\documentclass[11pt,a4paper]{amsart}
\usepackage{mathrsfs,bbm}
\usepackage[all]{xy}

 \newtheorem{thm}{Theorem}[section]
 \newtheorem{cor}[thm]{Corollary}
 \newtheorem{lem}[thm]{Lemma}
 
 \theoremstyle{definition}
 
 \theoremstyle{remark}

 \numberwithin{equation}{section}

\def\nn{^{(m)}}		
\def\xx{\bar{x}}	
	\def\XX{\times}
	\def\mP{\mathbb P}
\def\cir{\mathop{\stackrel{\scriptscriptstyle\circ}{}}}
\def\gl{\mathop\mathrm{GL}\nolimits}
		
\def\mF{\mathbb F}	
	\def\8{\infty}
	\def\+{\oplus}
	\def\aK{\mathbbm k}
\def\la{\lambda}	\def\mN{\mathbb N}
\def\gG{\mathfrak g}	
	
		\def\2{^{(2)}}	
\def\dA{\mathfrak{A}}	\def\mzg{\mathbb{Z}_{\ge0}}
\def\ka{\kappa}		\def\dF{\mathfrak{F}}
\def\tP{\tilde{\mP}}	\def\eps{\varepsilon}
\def\tF{\tilde{\dF}}	\def\tA{\tilde{\dA}}
\def\dR{\mathfrak{R}}
\def\res{\mathop\mathrm{res}}
\def\co{\mathop\mathrm{co_1}}

\def\gnrsuch#1#2{\langle\,#1\mid #2\,\rangle}
\def\lst#1#2{ #1_1 , #1_2 , \dots , #1_{#2} }
\def\row#1#2{( #1_1 , #1_2 , \dots , #1_{#2} )}
\def\bop{\bigoplus}
\def\mtr#1{\begin{pmatrix}#1\end{pmatrix}}

\def\FF{\phantom{\text{\LARGE I}}}
\def\gl{\mathop\mathrm{GL}}

\def\chg{Cher\-ni\-kov $p$-group}
\def\cht{Cher\-ni\-kov $2$-group}
\def\iff{if and only if }

\def\qA{{\boldsymbol A}}		\def\qR{{\boldsymbol R}}
\def\qB{{\boldsymbol B}}		\def\qE{\boldsymbol E}
				\def\mF{\mathbb F}
						
\def\aK{\mathbbm k}

\def\+{\oplus}		\def\xx{\times}

\def\set#1{\left\{#1\right\}}
\def\setsuch#1#2{\{\,#1\mid #2\,\}}

\begin{document}
\title[Nilpotent Chernikov $2$-groups]{On nilpotent Chernikov $2$-groups with elementary tops}
\author[Y. Drozd]{Yuriy A. Drozd}
\author[A. Plakosh]{Andriana I. Plakosh}

\address{Institute of Mathematics, National Academy of Sciences of Ukraine, 01601 Kyiv, Ukraine}
\email{y.a.drozd@gmail.com,\,drozd@imath.kiev.ua}
\email{andrianalomaga@mail.ru}
\urladdr{http://www.imath.kiev.ua/$\sim$drozd}

\maketitle

\section{Introduction}
\label{s1}

 Recall that a \emph{\chg} \cite{ch,kur} $G$ is an extension of a finite direct sum $M$ of \emph{quasi-cyclic} $p$-groups, 
 or, the same, the groups of type $p^\8$, by a finite $p$-group $H$. Note that $M$ is the biggest abelian 
 divisible subgroup of $G$, so both $M$ and $H$ are defined by $G$ up to isomorphism. We call $H$ and $M$,
 respectively, the \emph{top} and the \emph{bottom} of $G$. We denote by $M\nn$ a direct sum
 of $m$ copies $M_k\ (1\le k\le m)$ of quasi-cyclic $p$-groups and fix elements $a_k\in M_k$ of order $p$. 
 The group $G$ is nilpotent \iff the induced action of $H$ on $M$ is trivial \cite[Theorem~1.9]{ch}.
 
 In the papers \cite{sha,dp} the classification of nilpotent \chg{s} with elementary tops was related to 
 the classification of tuples of skew-symmatric matrices over the filed $\mF_p$. Namely, given an $m$-tuple of $n\xx n$ 
 skew-symmetric matrices $\qA=(\lst Am)$, where $A_k=(a^{(k)}_{ij})$, we define the \chg\ $G(\qA)$, which is an 
 extension of $M\nn$ by the elementary $p$-group $H_n=\gnrsuch{\lst hn}{h_i^p=1,\,h_ih_j=h_jh_i}$ such that $[h_i,a]=1$ 
 for each $a\in M\nn$ and $[h_i,h_j]=\sum_ka^{(k)}_{ij}a_k$. Every nilpotent \chg\ is of this kind and two $m$-tuples 
 $\qA=\row Am$ and $\qB=\row Bm$ define isomorphic groups \iff there are invertible matrices $S\in\gl(n,\mF_p)$ and
 $Q=(q_{kl})\in\gl(m,\mF_p)$ such that $B_k=\sum_l q_{lk}(SA_lS^\top)$ for all $k$. In this case we write $\qB=S\cir\qA\cir Q$
 and call the $m$-tuples $\qA$ and $\qB$ \emph{weakly equivalent}. Recall that the pairs $\qA$ and $S\cir\qA$ are called 
  \emph{congruent}.
 
 If $m>2$, a classification of $m$-tuples of skew-symmetric matrices is a \emph{wild problem} in the sense of the representation
 theory, i.e. it contains a classification of representations of any finitely generated algebra\cite{dp}. So, there is no hope to obtain a ``good''
 classification of \chg{s} with the bottom $M\nn$ for $m>2$.
 Using the results of \cite{ser}, we gave in the paper \cite{dp} a classification of \chg{s} with elementary tops and
 the bottom $M^{(2)}$ for $p\ne2$. Unfortunately, if $p=2$, the technique of \cite{ser} does not work. In this paper we use instead
 the results of \cite{wat} to obtain an analogous classification for \cht{s}.
 
\section{Alternating pairs}
 \label{s2} 

From now on $\aK$ is a field of characteristic $2$. We consider pairs $(A, B)$
of alternating bilinear forms in a finite dimensional vector space over $\aK$
or, the same, pairs of skew-symmetric matrices over $\aK$, calling them
\emph{alternating pairs}. Let $\qR = \aK[t]$, the polynomial ring, $\qE = \aK(t)/\aK[t]$ and
$\res = \res_{\infty} : \qE \to \aK$ be the residue at infinity. Let $M$ be a finite dimensional
(over $\aK$) $\qR$-module and $F : M \times M  \to \qE$ be an $\qR$-bilinear map. We call
$F$ \emph{strongly alternating} if $\res F(u, u) = \res F(tu, u) = 0$ for all $u \in M$. Then
also $F(u, v) = F(v, u)$ and $F(tu, v) = F(tv, u)$. Given a strongly alternating
map $F$ we set $A_F(u, v) = \res F(u, v)$ and $B_F(u, v) = \res F(tu, v)$. Obviously,
$(A_F, B_F)$ is a pair of alternating bilinear forms on $M$. We use the following
facts from \cite{wat}.

\textbf{Fact 1.} The map $F \mapsto (A_F, B_F) $ induces a one-to-one correspondence between
isomorphism classes of non-degenerated strongly alternating maps and isomorphism
classes of pairs of alternating forms $(A, B)$ such that $A$ is non-degenrated.

\textbf{Fact 2.} Isomorphism classes of indecomposable non-degenerated strongly alternating
maps $F : M\times M \to \qE$ are in one-to-one correspondence with powers
$f^n(t)$ of irreducible polynomials $f(t) \in \aK[t]$. Namely $f^n(t)$ corresponds to
the strongly alternating map $F_{f,n} : M_{f,n}\to \qE$, where $M_{f,n} = {(\qR/f^n\qR)}^2=
\gnrsuch{u, v}{f^nu=f^nv=0}$, such that $F_{f,n}(u, v) = 1/f^n(mod \aK[t])$, while
$F_{f,n}(u, u) = F_{f,n}(v, v) = 0.$

We denote the alternating pair corresponding to the map $F_{f,n}$ by $\qA_{f,n}=(A_{f,n}, B_{f,n}).$
Consider the matrices of size $n \times (n + 1)$
$$I_{n+}=\left(
       \begin{array}{ccccc}
         1 & 0 & 0 & \ldots & 0 \\
         0 & 1 & 0 & \ldots & 0 \\
         \hdotsfor{5} \\
         0 & 0 & 0 & \ldots & 1 \\
         0 & 0 & 0 & \ldots & 0 \\
       \end{array}
     \right)\!,\quad 
  I_{n-}=   \left(
       \begin{array}{cccccc}
         0 & 0 & 0 & \ldots & 0 &0\\
         1 & 0 & 0 & \ldots & 0 &0\\
         0 & 1 & 0 & \ldots & 0 &0\\
         \hdotsfor{6} \\
         0&0&0& \dots & 1 &0 \\
         0 & 0 & 0 & \ldots &0& 1 \\
       \end{array}
       \right)$$
and alternating pairs 
\begin{align*}
 \qA_{\infty,n}&=(A_{\infty,n}, B_{\infty,n}),\quad  \qA_{+,n}=(A_{+,n}, B_{+,n}), \notag\\
 \intertext{where}
A_{\infty,n}&=\left(
       \begin{array}{cc}
         0 & J_n  \\
         J_n^\top & 0  \\
       \end{array}
       \right),\quad
B_{\infty,n}=\left(
       \begin{array}{cc}
         0 & I_n  \\
         I_n & 0  \\
       \end{array}
       \right) \\
       A_{+,n}&=\left(
       \begin{array}{cc}
         0 & I_{n+}  \\
         I^\top_{n+} & 0  \\
       \end{array}
       \right),\quad
B_{+,n}=\left(
       \begin{array}{cc}
         0 & I_{n-}  \\
         I^\top_{n-} & 0  \\
       \end{array}
       \right), \notag
\end{align*}
$I_n$ is the $n\times n$ unit matrix and $J_n$ is the $n\times n$ nilpotent Jordan block.

\textbf{Fact 3.} Every indecomposable alternating pair $(A, B)$ with the degenerated form
$A$ is isomorphic to one of the pairs $(A_{\infty,n}, B_{\infty,n}), (A_{+,n}, B_{+,n}),$.

\textbf{Fact 4.} Every alternating pair decomposes into an orthogonal direct sum of
indecomposable pairs. This decomposition is unique up to isomorphism and
permutation of summands.

 \begin{lem} \label{11}
There is a $\aK$-basis in $M_{f,n}$ such that the forms $A_{f,n}$ and $B_{f,n}$
are given by the matrices $A_{f,n}=\left(
       \begin{array}{cc}
         0 & I  \\
         I & 0  \\
       \end{array}
       \right)$  a
and $B_{f,n}=\left(
       \begin{array}{cc}
         0 & \Phi  \\
         \Phi^\top & 0  \\
       \end{array}
       \right),$ where $\Phi$
is the Frobenius matrix with the characteristical polynomial $f^n(t)$.
 \end{lem}

Note that  $(A_{\infty,n}, B_{\infty,n})=(B_{t,n},A_{t,n})$.

\begin{proof}
We include $\aK[t]$ into the ring $\aK[[t]]$ of formal power series and into
the field $\aK((t))$ of Laurent series. If $\deg g = d$ and $g(0)\neq 0$, we set $g^*(t)=t^dg(1/t)$
and choose a polynomial $\tilde{g}(t)$ of degree $d$ such that $g^*(t)\tilde{g}(t)\equiv1\pmod{ t^{d+1}}$.
 It exists and is unique since $g^*(t)$ is invertible in $\aK[[t]]$.

Let $f(t)\neq t,\, g(t)=f^n(t),\, d=\deg g(t)$ and $g(t)=t^d+\alpha_1 t^{d-1}+\ldots +\alpha_d.$
Then $g^*(t)=1+\alpha_1t+\ldots+\alpha_d t^d$ and $\tilde{g}(t)=1+\beta_1t+\ldots+\beta_d t^d, $ 
where, for every $m\leq d,$
\begin{equation}\label{e1}
\alpha_m+\alpha_{m-1}\beta_1+\alpha_{m-2}\beta_2+\ldots+\alpha_1\beta_{m-1}+\beta_m=0
\end{equation}
(we set $\alpha_0=\beta_0=1$). Consider the basis $\setsuch{ u_k,v_k }{ 0\leq k < d}$ of $M_{f,n}$, 
where $v_k=t^kv, u_k=t^{d-k-1}u.$
Then $F_{f,n}(u_k,u_l)=F_{f,n}(v_k,v_l)=0$ for all $k,l,$ 
while $F_{f,n}(u_l,v_k)=h_{k,l}=t^{d+k-l-1}/ g(t)\pmod{ \aK[[t]]}. $
 Denote by $\co h$ the coefficient by $t^{-1}$ in the Laurent series $h$.
Recall that $\res_{\infty}h$, where $h \in \aK((t))$, equals $\co t^{-2}h(1/t)$. Therefore,
\begin{align*}
A_{f,n}&=\co t^{-2}h_{k,l}(1/t)=\co \frac{t^{l-k-1}}{t^d g(1/t)}=\co t^{l-k-1}\tilde{g}(t)=\begin{cases}
{\beta}_{k-l} & \text{if $k \geq l$,} \\
0 & \text{if $k < l$;}
\end{cases}\\
B_{f,n}&=\co t^{-3}h_{k,l}(1/t)=\co \frac{t^{l-k-2}}{t^d g(1/t)}=\co t^{l-k-2}\tilde{g}(t)=\begin{cases}
{\beta}_{k-l+1} & \text{if $k \geq l-1$,} \\
0 & \text{if $k < l-1$.}
\end{cases}
\end{align*}
So the matrices of the forms $A_{f,n}$ and $B_{f,n}$ in this basis are, respectively,
\begin{equation}\label{e2}
\left(
       \begin{array}{cc}
         0 & A  \\
         A^\top & 0  \\
       \end{array}
       \right)
 \text{and}      
\left(
       \begin{array}{cc}
         0 & B  \\
         B^\top & 0  \\
       \end{array}
       \right),       
\end{equation}
where
\begin{align*}
A&=\left(
       \begin{array}{ccccc}
         1 & \beta_1 & \beta_2 & \ldots & \beta_{d-1} \\
         0 & 1 & \beta_1 & \ldots & \beta_{d-2} \\
         0 & 0 & 1 & \ldots & \beta_{d-3} \\
          \hdotsfor{5} \\
         0 & 0 & 0 & \ldots & 1 \\
       \end{array}
     \right),\\
B&=\left(
       \begin{array}{cccccc}
          \beta_1 & \beta_2 & \beta_3 & \ldots & \beta_{d-1} & \beta_d \\
          1 & \beta_1 & \beta_2 & \ldots & \beta_{d-2} & \beta_{d-1} \\
          0 & 1& \beta_1 & \ldots & \beta_{d-3} & \beta_{d-2}\\
          \hdotsfor{6} \\
         0 & 0 & 0 & \ldots & 1 & \beta_1 \\
       \end{array}
     \right).\\
\intertext{The relations (\ref{e1}) imply that}
     A^{-1}&=\left(
       \begin{array}{ccccc}
         1 & \alpha_1 & \alpha_2 & \ldots & \alpha_{d-1} \\
         0 & 1 & \alpha_1 & \ldots & \alpha_{d-2} \\
         0 & 0 & 1 & \ldots & \alpha_{d-3} \\
          \hdotsfor{5} \\
         0 & 0 & 0 & \ldots & 1 \\
       \end{array}
     \right),
\end{align*}
     and $A^{-1}B=\Phi$, the Frobenius matrix with the characteristical polynomial 
     $g(t)=f^n(t).$ Thus, multiplying the matrices of bilinear forms $A_{f,n}$ and $B_{f,n}$ from (\ref{e2})
     by the matrix
     $$\left(
       \begin{array}{cc}
         A^{-1} & 0 \\
         0 & I  \\
       \end{array}
       \right)
     $$
     on the left and by the transposed matrix on the right, we accomplish the
proof of the lemma in this case.

If $f(t) = t$, we obtain the necessary form of the matrices directly in the
basis $\{u_k, v_k\}$ as above.
\end{proof}

Now we resume the above considerations.

\begin{thm} 
Every indecomposable alternating pair is isomorphic to one
of the pairs
\[
 \qA_{f,n}=(A_{f,n},B_{f,n}),\, \qA_{\infty,n}=(A_{\infty,n},B_{\infty,n}),\,\qA_{+,n}= (A_{+,n},B_{+,n})
 \]
 given by Fact~3 and Lemma~\ref{11}.
Every alternating pair decomposes uniquely (up to permutation of summands)
into an orthogonal sum of indecomposable strongly alternating pairs
from this list.

\end{thm}

\section{Weak equivalence and Chernikov groups}
\label{s3} 

We denote by $\dA$ the set of all pairs $\qA$, where 
$\qA\in\set{\qA_{f,n},\qA_{\infty,n},\qA_{+,n}}$, and by $\dF$
 the set of functions $\ka:\dA\to\mzg$ such that $\ka(\qA)=0$ for almost all $\qA$. For any function
 $\ka\in\dF$ we set $\dA^\ka=\bop_{\qA\in\dA}\qA^{\ka(\qA)}.$ For the classification of \cht{s} we have to answer
 the question:
 
 \smallskip
 \emph{Given two functions  with finite supports $\ka,\ka':\dA\to\mzg$, when are the pairs $\dA^\ka$ and $\dA^{\ka'}$
 weakly congruent}?

 \smallskip
 Evidently, $(\qA_1\+\qA_2)\cir Q=(\qA_1\cir Q)\+(\qA_2\cir Q)$, so the pairs $\qA$ and $\qA\cir Q$ are indecomposable
 simultaneously. For every pair $\qA\in\dA$ we denote by $\qA*Q$ the unique pair from $\dA$ which is congruent to $\qA\cir Q$.
 The map $\qA\mapsto\qA*Q$ defines an action of the group $\gG=\gl(2,\aK)$ on the set $\dA$, hence on the set $\dF$ of functions
 $\ka:\dA\to\mzg$: $(Q*\ka)(\qA)=\ka(\qA*Q)$.

 \begin{cor}\label{21}
  The pairs $\dA^\ka$ and $\dA^{\ka'}$ are weakly congruent \iff the functions $\ka$ and $\ka'$ belong to the same
  orbit of the group $\gG$.
 \end{cor}

$ (A_{+,n},B_{+,n})$ is a single indecomposable couple of dimension $2n-1$.
For every other pair $\qA=(A,B)$ the polynomial $\det (xA+yB)$ is a square:
$\det (x_1A+x_2B)=\Delta_{\qA} (x_1,x_2)^2$ for some $\Delta_{\qA}(x_1,x_2)$ (the \emph{Pfaffian} of $x_1A+x_2B$, 
see \cite{km}). Namely, 
\[
 \Delta_\qA(x,y)=
 \begin{cases}
  x_2^n &\text{ if } \qA=\qA_{\8,n},\\
  x_2^{dn}f(x_1/x_2) &\text{ if } \qA=\qA_{f,n} \text{ and } \deg f=d.
 \end{cases}
\]
If $(A',B')=(A,B)\circ Q$, where $Q=\mtr{q_{11} & q_{12} \\ q_{21} & q_{22} }$, 
       then $\Delta_{(A',B')}(x_1,x_2)=\Delta_{(A,B)}((x_1,x_2)Q)=\Delta_{(A,B)}(q_{11}x_1+q_{21}x_2,q_{12}x_1+q_{22}x_2)$.
       So now we can repeat the considerations of \cite{dp}, obtaining analogous results for the fields of characteristic $2$
       and \cht{s}.
       
 We say that an irreducible homogeneous polynomial $g\in\aK[x_1,x_2]$ is \emph{unital} if either $g=x_2$ or its leading
 coefficient with respect to $x_1$ equals $1$. Let $\mP=\mP(\aK)$ be the set of unital homogeneous irreducible polynomials
 from $\aK[x_1,x_2]$ and $\tP=\tP(\aK)=\mP\cup\{\eps\}$. Note that $\mP$ actually coincides with the set of the closed
 points of the projective line $\mP^1_{\aK}=\mathrm{Proj}\,\aK[x_1,x_2]$ \cite{ha}. For $g\in\mP$ and $Q\in \gG$,
 let $Q*g$ be the unique polynomial $g'\in\mP$ such that $g((x,y)Q)=\la g'$ for some non-zero $\la\in\aK$.
 (It is the natural action of $\gG$ on $\mP^1_{\aK}$.) We also set $Q*\eps=\eps$ for any $Q$. It defines an action
 of $\gG$ on $\tP$. Denote by $\tF=\tF(\aK)$ the set of all functions $\rho:\tP\XX\mN\to\mzg$ such that $\rho(g,n)=0$
 for almost all pairs $(g,n)$. Define the actions of the group $\gG$ on $\tF$ setting $(\rho*Q)(g,n)=\rho(Q*g,n)$.
 For every pair $(g,n)\in\tP\xx\mN$ we define a pair of skew-symmetric forms $\qA(g,n)$:
 \[
  \qA(g,n)=\begin{cases}
  (A_{\infty,n}, B_{\infty,n}) &\text{if } g=x_2,\\
  (A_{+,n}, B_{+,n}) &\text{if } g=\varepsilon,\\
    (A_{f,n},B_{f,n})&\text{where } f=g(x,1) \text{ otherwise}.
  \end{cases}
 \]
 Let $\tA=\tA(\aK)=\{\qA(g,n)\,|\,(g,n)\in\tP\XX\mN\}$. For every function $\rho\in\tF$ we set
 $\tA^\rho=\bop_{(g,n)\in\tP\XX\mN}\qA(g,n)^{\rho(g,n)}$.
 The preceding considerations imply the following theorem.

 \begin{thm}\label{22}
 \begin{enumerate}
 \item  Every pair of skew-symmetric bilinear forms over the field $\aK$ is weakly congruent to $\tA^\rho$ for some function
 $\rho\in\tF(\aK)$.
 \item  The pairs $\tA^\rho$ and $\tA^{\rho'}$ are weakly congruent \iff the functions $\rho$ and $\rho'$ belong to
 the same orbit of the group $\gG=\gl(2,\aK)$.
 \end{enumerate}
 \end{thm}

For every function $\rho\in\tF(\mF_2)$ set $G(\rho)=G(\tA^\rho)$.

 \begin{thm}\label{23}
 Let $\dR$ be a set of representatives of orbits of the group $\gG=\gl(2,\mF_2)$ acting on the set of functions
 $\tF(\mF_p)$. Then every nilpotent \cht\ with elementary top and the bottom $M\2$ is isomorphic to the group
 $G(\rho)$ for a uniquely defined function $\rho\in\dR$.
 \end{thm}

The description of these groups in terms of generators and relations is also the same as in \cite{dp}. Note that all of them are of the form
 $G(\qA)$, where $\qA=\bop_{k=1}^s\qA_k$ and all $\qA_k$ belong to the set $\set{\qA_{\infty,n},\qA_{+,n},\qA_{f,n}}$. Each term $\qA_k$
 corresponds to a subset $\set{h_{ki}}$ of generators of the group $H$ and we have to precise the values of $[h_{ki},h_{kj}]$
 (all other commutators are zero). They are given in Table~1. Recall that $a_1$ and $a_2$ are generators of the subgroup
 $\setsuch{a\in M^{(2)}}{2a=0}$.

 \begin{table}[!ht]
\caption{}\vspace*{-1em}
 \[
  \begin{array}{|c|c|c|}
  \hline
  \FF \qA_k & i,j & [h_{ki},h_{kj}] \\
  \hline
    \FF \qA_{+,n}\ & j=d+i & a_1\\
  & j=d+i-1 & a_2 \\
  & \text{ otherwise } & 0  \\
  \hline
  \FF \qA_{\infty,n} & j=d+i & a_2,\\
  & j=d+i-1 & a_1,\\
  & \text{ otherwise } & 0 \\
  \hline
  \FF \qA_{f,n} &\ j=d+i<2d\ & a_1\\
   &j=d+i-1 &
  a_2 \\
  & i<d,\,j=2d &\  \la_{d-i+1}a_2\ \\
  & i=d,\,j=2d &\ a_1+\la_1a_2\  \\
  & \text{ otherwise } & 0  \\
  \hline
  \end{array}
 \]
 
 \smallskip
 \centerline{ where $f^n(x)=x^d+\la_1 x^{d-1}+\dots+\la_d$}
 \end{table}

\newpage

\label{lit}

\end{document}